\documentclass[11pt,draft]{article}
\usepackage[final]{graphics}
\usepackage{hyperref}
\usepackage{graphics}
\usepackage{amsfonts}
\usepackage{amsmath,amsfonts,mathrsfs,amssymb,color}
\usepackage{indentfirst}
\usepackage{graphicx}
\usepackage{float}

\numberwithin{equation}{section}

\newtheorem{theo.}{\quad\, Theorem}[section]
\newtheorem{defi.}{\quad\, Definition}[section]
\newtheorem{lemm.}{\quad\, Lemma}[section]
\newtheorem{coro.}{\quad\, Corollary}[section]

\begin{document}

\title{A Mathematic Expression of the Genes of Chinese Traditional Philosophy$^*$ }
\author{Kegong Chen$^{a}$
\ \ \ \ \ \ \ \ \ Ruyun Ma$^{b}$
\\
 {
 \small\it $^{a,b}$Northwest
Normal University, Lanzhou 730070, P R China}
}
\date{} \maketitle
\footnote[0]{E-mail address: chenkg@nwnu.edu.cn (K. Chen);
mary@nwnu.edu.cn (R. Ma)
%\ \ mary@nwnu.edu.cn (R. Ma), \ \
%Telephone: 86-931-7971297
% $^{\ddagger} $xuling$_{-}$216@yahoo.cn
} \footnote[0] {$^*$Supported by the NSFC (No.11671322). }
 \begin{abstract}
We provide a mathematic model for the Traditional Yin-and-Yang Double Fish Diagram which from Chinese Traditional Philosophy.

\end{abstract}
\vskip 3mm

{\small\bf Keywords.} {\small Mathematic model,
Traditional Yin-and-Yang Double Fish Diagram, Chinese traditional philosophy, Pythagorean theorem}

\vskip 3mm

%{\small\bf MR(2000)\ \ \ 34B10, \ 34B18}

\baselineskip 24pt

\section{Introduction}

\vskip 3mm

The Book of Changes (also known as I Ching), profound and wide-spreading, is a classic work with a long history to represent the Chinese classic traditional philosophical thinking. In Yi Zhuan (Commentaries on the Book of Changes), (also known as Ten Wings, which was said to be written by Confucius), Confucius wrote ten articles to make the Book of Changes complete and more attractive. The Traditional Yin-and-Yang Double Fish Diagram, as a development of the ideas in I Ching, endows the traditional philosophy ``Yin and Yang give birth the ultimate truths" with a vivid expression, which, however, is quite different from Western logics. The difference makes the world think that the most essential Chinese traditional philosophy lack of rationale and scientific explanation. In this perspective, the Traditional Yin-and-Yang Double Fish Diagram has been regarded inexplicable and inexpressible because it has never been expressed and proved sufficiently by nay functional equation and has not been established as a standard diagram, which consequently posed a negative influence on the cultural confidence of the Chinese traditional culture.

The present paper starts with the Pythagorean Theorem. By way of the inscribed circle within unit squares and the application of the core ideas of the Lever Balance Principle (Archimedes Law) and Fuzzy Mathematics, the writers successfully express the $S$-shaped curve of the Traditional Yin-and-Yang Double Fish Diagram with an elementary function, thus successfully explain the unspeakable meaning of the question in mathematic language. This study ends the history that there was no standard equation and no standard geometrical figure to illustrate the Traditional Yin-and-Yang Double Fish Diagram. The significance of the present study lies not only in sufficiently expressing the scientific value of the Traditional Yin-and-Yang Double Fish Diagram, but also in explaining the fundamental basis on which the Traditional Yin-and-Yang Double Fish Diagram reflects the nature of all existence and development.

\section{The Philosophical Content of the Pythagorean Theorem}

The usual definition for a circle is: when a line segment rotates around one of its two end points on a plane, the trajectory of the other end point is called a CIRCLE. Another definition goes like: in the same plane, the set of points which have equal distance to a fixed point compose a CIRCLE. The definition of a circle can also be obtained from the well-known Pythagorean Theorem---In a right angled triangle, the square of the hypotenuse is equal to the sum of the squares of the other two sides. In other words, when the length of the hypotenuse is fixed, the trajectory of the end points of the other two sides is called a CIRCLE. This is called the Gou Gu Theorem in China while the Pythagorean Theorem in the Western world. The following Geometrical Figure 1 and Algebraic Equation (1) are illustrations of the theorem:

%\begin{center}
%\includegraphics[height=40mm]{image/2.pdf}
%\end{center}
\vskip 20mm
\begin{figure}[!h]
\begin{center}
\scalebox{0.4}{\includegraphics{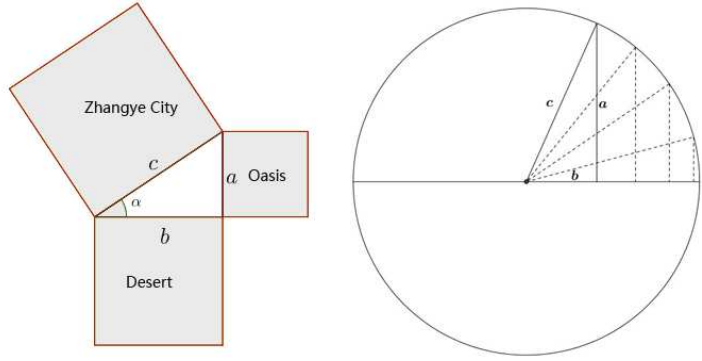}}
\end{center}
\caption{The Pythagorean Theorem and Geometrical Figure of Circle.}
\label{fig1}
\end{figure}

Figure 1 means:
$$
a^2+b^2=c^2,\ \ \ \text{or}\ \ \ \ \Big(\frac{a}{c}\Big)^2+\Big(\frac{b}{c}\Big)^2=1.
\eqno(1)$$

Zhao Shuang, a mathematician in the late Eastern Han Dynasty, remarked when he annotated Zhou Bi Suan Jing, an ancient Chinese classic on mathematics, that ``The sum of the areas of the two squares on the legs ($a$ and $b$) equals the area of the square on the hypotenuse ($c$).''

We use $x$ and $y$ to represent the quotient of the square of two right-angle sides and the square of the hypotenuse respectively, and define by $x$ and $y$ the value of Yang and the value of Yin of the right angled triangle respectively:
$$x=\frac{a^2}{c^2},\ \ \ y=\frac{b^2}{c^2}.$$
Then
$$x+y=1.\eqno(2)$$

This means however the two right angle sides vary, the sum of the quotient values is identically equal to 1. That is, the sum of the quotient of the value of Yang and the quotient of the value of Yin is identically 1. Its geometric form is as Figure 2:

\vskip 20mm
\begin{figure}
\begin{center}
\scalebox{0.4}{\includegraphics{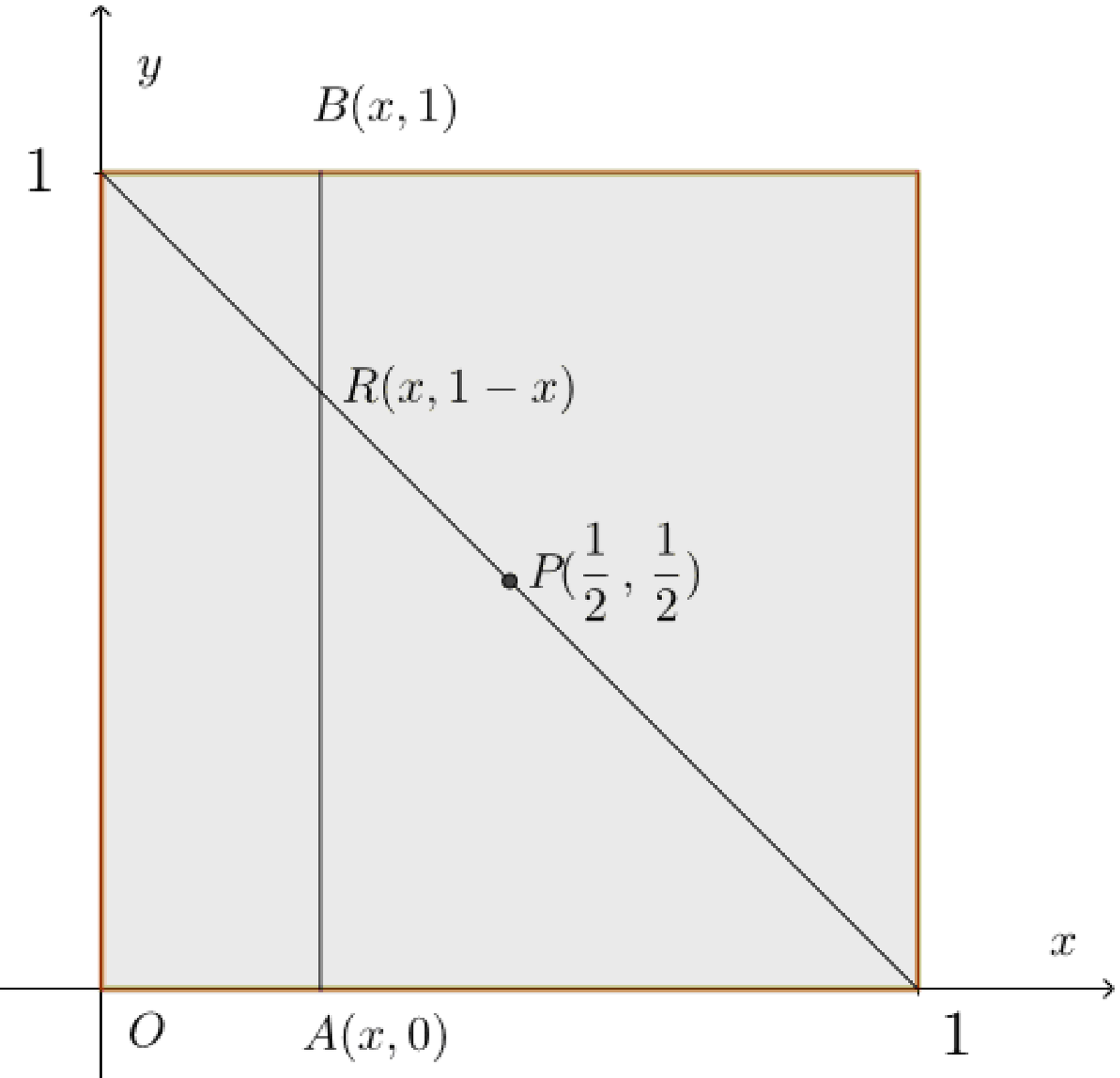}}
\end{center}
\label{fig1}
\caption{}
\end{figure}

Since\ $a, b\in(0, c)$, it follows that for $x, y\in(0, 1)$,
$$x=\frac{a^2}{c^2}=\frac{c^2-b^2}{c^2}=\sin^2\alpha=\frac{1-\cos2\alpha}{2}:=P(\alpha),\ \
\alpha\in(-\infty,\ +\infty),$$
$$x=\frac{b^2}{c^2}=\frac{c^2-a^2}{c^2}=\cos^2\alpha=\frac{1+\cos2\alpha}{2}:=N(\alpha),\ \
\alpha\in(-\infty,\ +\infty),\eqno(3)$$

\begin{figure}[!h]
\begin{center}
\scalebox{0.4}{\includegraphics{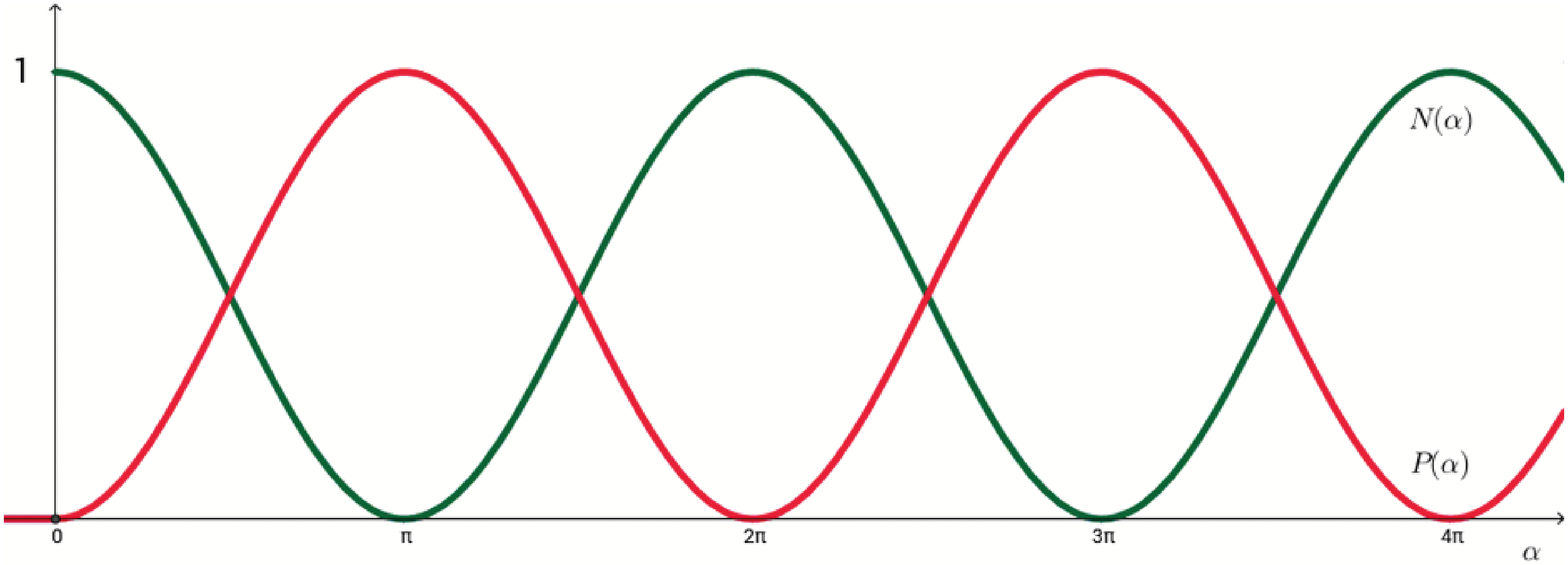}}
\end{center}
\caption{}
\end{figure}

Algebraic Equation (1) and Figure 1 are the foundation of Euclidean Geometry, while Figure 2 and Algebraic Equation (2) suggested the idea of Fuzzy Set, which paved the way for the development of modern mathematics. Combining them, the concept can be reached as follows:

The essence of the Pythagorean Theorem is that the sum of the areas of the two squares on the legs equals the area of the square on the hypotenuse. We can elaborate on this in relation to the squares of the oasis and desert in Zhangye City, Gansu, China. Suppose the total area of Zhangye City is 1, then the sum of the square of the oasis and the square of the desert is identically 1 (Zhangye City is so very unique that the oasis refers to the area with water, while the desert refers to the area without water. Besides, the sum of the oasis area and the desert area doesn't change). When the oasis is expanding, the desert will be shrinking, and vice versa. The shift in the set point is just the path of circle.

Since  $x=1-y$,  $y=1-x.$
Then the function of Yang and the function of Yin are both opposite and complementary, and a balanced unity of opposites. That is, each is based on the existence of the other party, while its total value is 1. In other words, the area of the oasis is the deviation of the area of the desert, and the area of the desert is the deviation of the area of the oasis. The sum of both deviations is identically 1. However, most areas in Zhangye City are neither absolutely desert nor oasis. Instead, most areas are between the oasis and the desert. It is just like the points between the interval $[0, 1]$ of the line segment $x+y=1$. What is different is that where the point is located. That is, whether it is closer to the left or to the right. Or, whether it is closer to the oasis or the desert. This is what the present author repeatedly highlighted in previous two articles (see more in Natural Science Edition of the Journal of Northwest Normal University, Issue 5 and Issue 6, 2015). For a system, we can set and classify the complementary opposites in a simple attribute, as the soil moisture can be divided into two sets and the state of the economy can be divided into two sets of consumption and accumulation. Greatness is in simplicity. In its simplest sense, we can treat all things quite subjectively, and according to objective equilibrium value, set and classify all the opposite and complementary things into positive and negative, white and black, up and down, left and right, big and small, and even good and bad (we can mark this classification method as  $\uparrow0\downarrow$). Even so, the real state of things is between the interval of white and black instead of being absolutely white or black (Russell's paradox). In most cases, it is both white and black. The only difference is that the proportion of white and black varies, which is also how the Fuzzy Mathematics developed.

Therefore, hypotenuse, the Pythagorean Theorem and the essence of circle are the same. They are philosophy and science about 1. Actually, the metaphysics that originated in the Wei and Jin Period is a philosophical thinking based on the Pre-Qin Yi thoughts, such as ``Yin and Yang are called Tao", ``Tao gives birth to One. One gives birth to Two. Two gives birth to Three. Three gives birth to all things." It is also based on the analytical thought of numbers and shapes relating to circle and the Pythagorean Theorem in Zhou Bi Suan Jing and Jiu Zhang Suan Shu, which deeply influenced the Chinese culture ever since. Chinese ink paintings are a typical representation of the thought. If black and white are defined respectively as 0 and 1, or Yin and Yang, then the universe will be represented between white and black, or in the interval $[0, 1]$. This vivid representation of gradual change brings people into the boundless imagination of the universe. Thus, we are confident that Chinese ink paintings are a map of Yin and Yang as well as a chart of Fuzzy Mathematics.

\section{The Internal Relations of the Pythagorean Theorem and the Traditional Yin-and-Yang Double Fish Diagram}
\vskip 3mm

In Figure 2, draw an inscribed circle of the square with $P$ as the center of the circle and $\frac{1}{2}$ length of a side as the radius, thus Figure 4.

\begin{figure}
\begin{center}
\scalebox{0.4}{\includegraphics{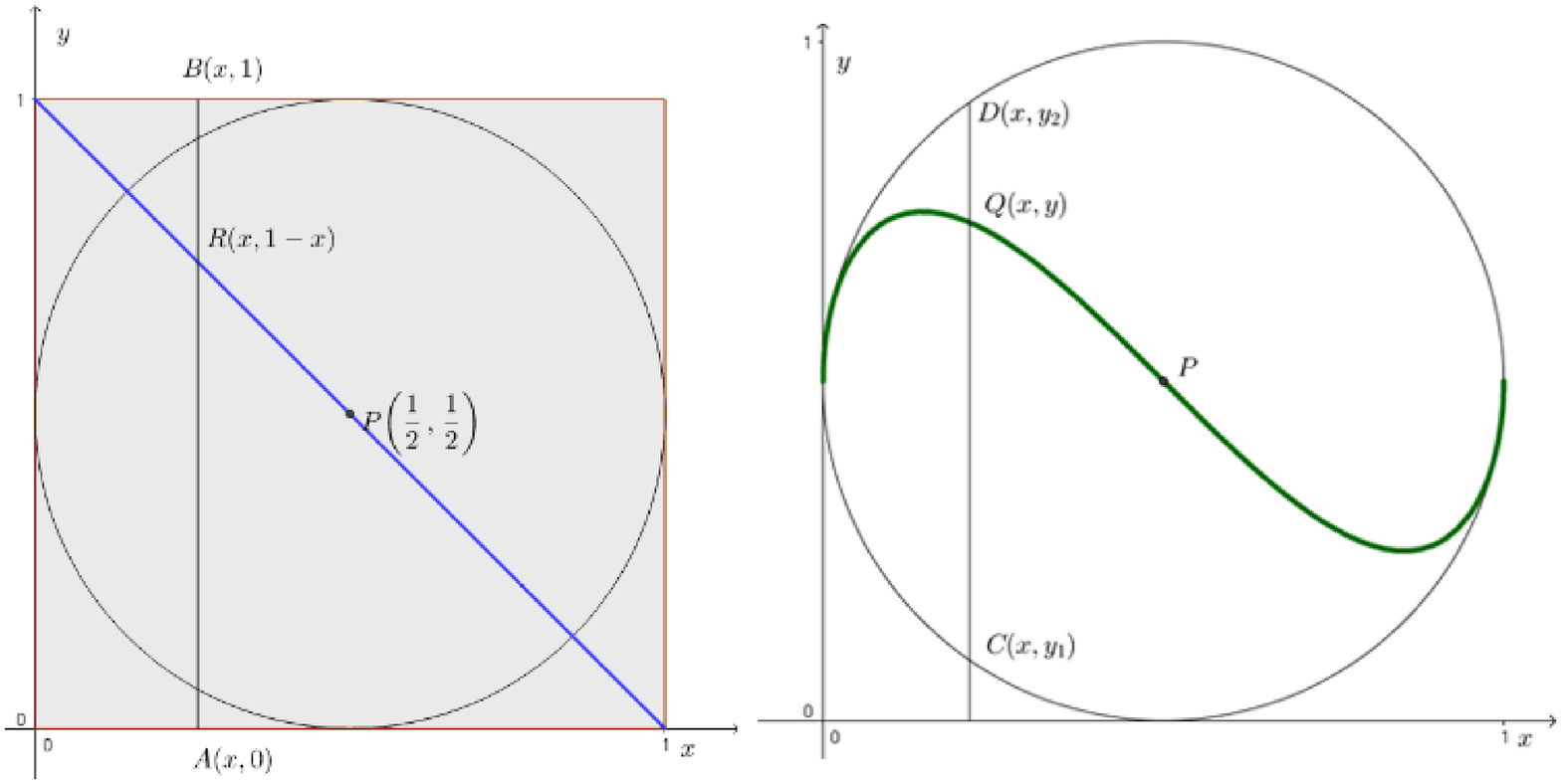}}
\end{center}
\caption{}
\label{fig1}
\end{figure}

In Figure 4, if $AB$ is considered as a level of standard unit, the length of $C(y_2-y)$ and
$D(y-y_1)$ are weights at each side, in the light of Lever Balance Principle, thus:
$$|AR|=1-x,\ \ \ |BR|=x,\ \ x\in[0, 1],$$
$$\frac{|AR|}{|BR|}=\frac{|CQ|}{|DQ|},\ \ \ \text{that\ is:}\ \
 \frac{1-x}{x}=\frac{y-y_1}{y_2-y}.$$
Then
$$y=xy_1+(1-x)y_2.\eqno(4)$$

If $C(x, y_1),\ D(x, y_2)$,\ then $C$   and $D$ satisfy the following circular equation:
$$\Big(x-\frac{1}{2}\Big)^2+\Big(y-\frac{1}{2}\Big)^2=\frac{1}{4},$$
In which $y_1=\frac{1}{2}-\sqrt{\frac{1}{4}-\Big(x-\frac{1}{2}\Big)^2},\
y_2=\frac{1}{2}+\sqrt{\frac{1}{4}-\Big(x-\frac{1}{2}\Big)^2},$\
substituting into Equation (4), thus we get the $S$-curve standard equation of the Traditional Yin-and-Yang Double Fish Diagram:
$$y=\frac{1}{2}+(1-2x)\sqrt{\frac{1}{4}-\Big(x-\frac{1}{2}\Big)^2},\ \ \ \
 x\in[0, 1].\eqno(5)$$
Take\ $y'=0$,\ thus $x_1=\frac{1}{2}-\frac{\sqrt{2}}{4},\
 x_2=\frac{1}{2}+\frac{\sqrt{2}}{4}$.

Thus $y$ has a maximum of
$\frac{3}{4}$ and a minimum of $\frac{1}{4}$. Consequently, we get the inflexion coordination of the Chinese Yin-and-Yang Diagram as  $\Big(\frac{1}{2}-\frac{\sqrt{2}}{4}, \frac{3}{4}\Big)$ and
$\Big(\frac{1}{2}+\frac{\sqrt{2}}{4}, \frac{1}{4}\Big)$ in Figure 4,
 and the ``fish eye" coordination are $\Big(\frac{1}{2}-\frac{\sqrt{2}}{4}, \frac{1}{2}\Big)$ and
$\Big(\frac{1}{2}+\frac{\sqrt{2}}{4}, \frac{1}{2}\Big)$.

Equation (4) illustrates the relations between the circle and the square. Although level and weights belong to different systems, those two systems are always keeping a renewed balance moving around a balance point with equal moment of force, i.e. although each side has different force arms and forces, the Yin and Yang rations are always equal.

 Equation (5) illustrates that when the line $x+y=1$ in system of rectangular coordinates is put in its inscribed circle, the line turns to a $S$-shape curve in Equation (5), as Figure 6:

\begin{figure}[!h]
\begin{center}
\scalebox{0.4}{\includegraphics{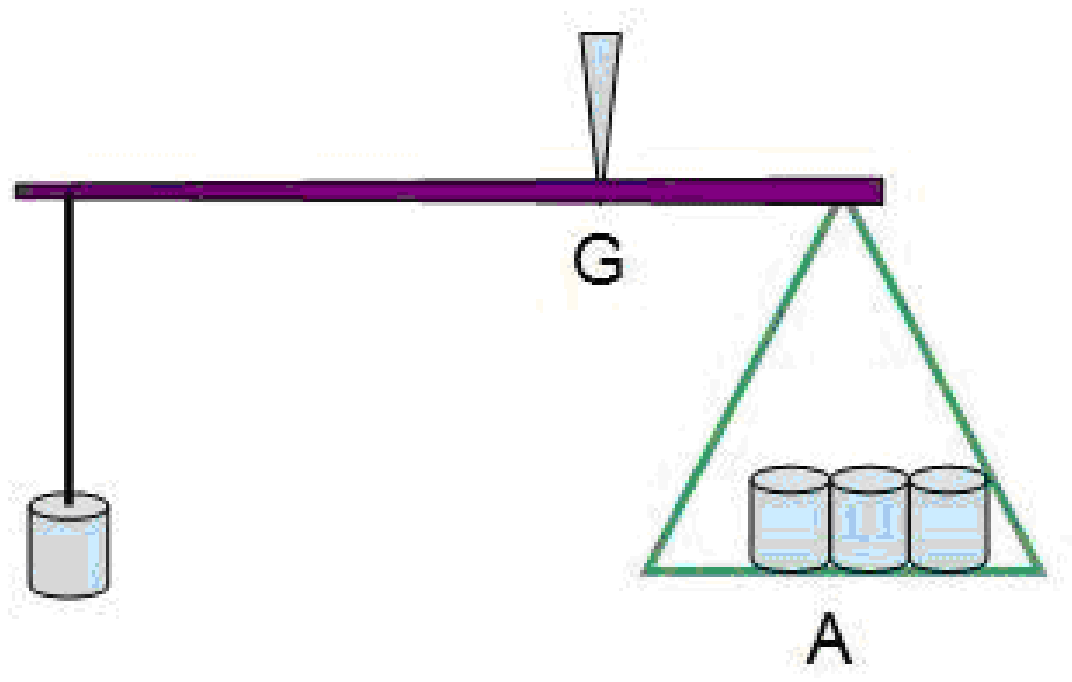}}
\end{center}
\caption{}
\label{fig1}
\end{figure}

\begin{figure}[!h]
\begin{center}
\scalebox{0.4}{\includegraphics{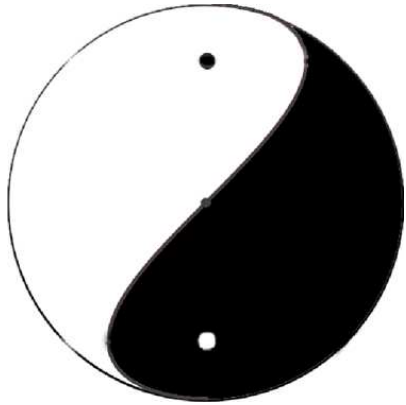}}
\end{center}
\caption{The Standard Traditional Yin-and-Yang Double Fish Diagram (The total of the area Yin and area Yang is 1, and the total of inverse values of Yin and Yang is 1)}
\label{fig1}
\end{figure}

\section{The Profound Meaning of the Traditional Yin-and-Yang Double Fish Diagram}

The universe is just like a right-angled triangle of standard unit. The two right-angle sides oppose each other and yet also complement each other like quantum entanglement. They superpose and counterbalance each other in a system. They grow and decline, but identically equal to ``1". From this perspective, the Traditional Yin-and-Yang Double Fish Diagram and the core idea behind it is counterbalance and unity of opposites, a geometric expression of ``being of beings" and an illustration of the ground of ``being", or the ``being" in ``being of beings", the ``one" in ``one is all". It roots in the Pythagorean Theorem and goes beyond the Pythagorean Theorem. It illustrates the idea in the Eular Equation: $e^{i\theta}\cdot e^{-i\theta}=(\cos\theta+i\sin\theta)(\cos\theta-i\sin\theta)=\cos^2\theta+\sin^2\theta=1$),
 and has opened the door of relativity, expressing the key idea of Fuzzy Set in modern mathematics developed in contemporary and modern times. ``If each right-angled side multiple itself, the sum will be the square of hypotenuse. The extracted root is hypotenuse". Thus, the xian (in Chinese meaning either ``hypotenuse" or ``profound") is ``profound" because we don't understand the principle and add too much mysterious explanations to it.

This reminds us of the ``Needham question", why didn't science rise in China? And the emotional sigh of Mr. Liang Shumin, ``The Chinese culture is a pre-matured one of mankind." In a material desire-pursuing physical world, metaphysics is so lonely and withering. As Heidegger pointed out, truth is hardly acceptable just because it's too simple. Facing the rapidly-changing advancement of science, we need inquire with earnestness and reflect with self-practice (quoted from the Analects), and ask ourselves whether we need ``lead people to perfection" and``having known where to reset at the end, one will be able to determine the object of pursuit" (quoted from the Great Learning), and try our best to ``return to things themselves"; Facing the cultural flourishing, whether we should keep to the original intention of ``speech" and come back to the original point of language.

Facing the complicated world, human beings need another renaissance to search, keep and return to the ``1" to create ``a community of common destiny". The philosophical element of traditional Chinese culture and the dialectical unity of Marxism could make continuing contributions to the aim since whether ancient or modern, Chinese or foreign, arts or sciences, they must be ``unified" in basic principle of the universe.

\vskip 10mm

\centerline {\bf REFERENCES}\vskip5mm\baselineskip 0.45cm
\begin{description}

\item{[1]} K. G. Chen,  R. Y. Ma, Mathematical model on Chinese Yin-and-Yang double-fish diagram and its application (I) 51(5) (2015) 1-3.　　
　

\item{[2]} K. G. Chen,  R. Y. Ma, Mathematical model on Chinese Yin-and-Yang double-fish diagram and its application (II) 51(6) (2015) 1-5.

\item{[3]} Jong-ping Hsu, Leonardo Hsu,　Space-Time Symmetry and Quantum Yang-Mills Gravity: How Space-Time Translational Gauge  Symmetry Enables the Unification of Gravity with Other Forces[M]. New Jersey: World Scientific, 2013.

\end{description}
\end{document}